\documentclass{conm-p-l}
\usepackage{comment,epsfig}

\newtheorem{theorem}{Theorem}[section]
\newtheorem{lemma}[theorem]{Lemma}
\newtheorem{corollary}[theorem]{Corollary}

\theoremstyle{definition}

\theoremstyle{remark}
\newtheorem{remark}[theorem]{Remark}

\numberwithin{equation}{section}

\def\SD{{\mathbb{S}}^{d-1}}
\def\SS{{\mathbb{S}}}
\def\RD{{\mathbb{R}}^{d}}
\def\RR{{\mathbb{R}}}
\def\NN{{\mathbb{N}}}

\def\eps{\epsilon}
\def\const{\qopname \relax o{const}}
\def\sgn{\qopname \relax o{sgn}}

\def\a{{\alpha}}
\def\ra{\rightarrow}

\def\NN{{\mathbb{N}}}

\def\min{\qopname \relax o{min}}

\def\pp{{p^\prime}}
\def\d{\delta}
\def\D{\Delta}

\begin{document}

\title[Moduli of continuity and average decay of Fourier transforms]
{Moduli of continuity and average decay of Fourier transforms:
two-sided estimates}

\author{Dimitri Gioev}
\address{Department of Mathematics, University of Rochester, Hylan Building,
Rochester, NY 14627}
\email{gioev@math.rochester.edu}
\thanks
{The author was supported in part by NSF grant INT--0204308
U.S.--Sweden Collaborative Workshop on PDE's and Spectral Theory,
the Swedish foundation STINT grant PD2001--128
and NSF grant DMS--0550649.
The author would like to thank the University of Pennsylvania
and the Courant Institute for financial support and hospitality.
The author is grateful to
Percy Deift for his valuable suggestions
which have helped to improve the presentation.
}
\date{December 15, 2006 and, in revised form, June 15, 2007.}


\begin{abstract}
We study inequalities between general integral moduli
of continuity of a function and the tail integral of its Fourier transform.
We obtain, in particular, a refinement of a result due 
to D. B. H. Cline \cite{Cl} (Theorem \ref{thm1} below).
We note that our approach does not use a regularly 
varying comparison function as in \cite{Cl}.
A corollary of Theorem 1.1
deals with the equivalence of
the two-sided estimates 
on the modulus of continuity on one hand,
and on the tail of the Fourier transform, on the other (Corollary
\ref{CH3_s3.lem.1}).
This corollary is applied in the proof of the violation of the so-called
entropic area law for a critical system of free fermions in \cite{GiKl,Gi2}.
\end{abstract}

\maketitle
\section{Introduction and statement of the main results}
\label{s0}
A result of this paper (Corollary~\ref{CH3_s3.lem.1} below) 
is applied in the proof of the violation of the so-called
entropic area law for a critical system of free fermions, 
see \cite[(6) et seq.]{GiKl}, \cite[Section on Fractal Boundaries]{GiKl}
 and \cite[Lemma 2.10]{Gi2}.
Corollary~\ref{CH3_s3.lem.1} follows from more
general results of this paper (Theorems \ref{thm1}, \ref{lem4.1})
which are of independent interest.

It is well-known that the behavior of a modulus of continuity
$\omega[f](h)$
of a function $f$ for $|h|$ small is related to the behavior of the 
Fourier transform $\hat{f}(\xi)$ of $f$ for $|\xi|$ large
(precise definitions are given in \eqref{eq200} et seq.~below),
see e.g. \cite[Proposition~5.3.4]{St}, \cite[Theorem~85]{T}.
The main object of our study are inequalities between 
general averaged moduli of continuity (\hbox{m.c.}) of $L^p$ functions
(defined in~\eqref{eq200} below)
and tails of their Fourier transforms (\hbox{F.t.}). 
In \cite{Cl} several results relevant for our purposes were obtained.
Theorem~\ref{thm1} below gives a lower estimate for a general $L^p$ 
\hbox{m.c.}, $1\leq{}p\leq2$, 
in terms of the modified tail integral of the \hbox{F.t.}
improving one of the results in~\cite{Cl} (as in \cite{Cl}, we
distinguish between the true and the modified
\hbox{F.t.} tail integral, as defined in \eqref{eq00true} and \eqref{eq00} below).
Corollary~\ref{thm2} 
gives a two-sided estimate for the \hbox{m.c.} in terms of
the modified tail integral of the \hbox{F.t.}
in the case $p=2$.
In applications it might be desirable to use the true \hbox{F.t.} tail
instead of the modified tail that  arises naturally in the mentioned
inequalities.
Theorem~\ref{lem4.1} gives the best possible
power-scale 
description of the 
relationship between the true and the modified \hbox{F.t.} tails
(see Remark~\ref{remstar}).

Before stating our results we need to introduce some notation
and recall two results in \cite{Cl}.
Let $d\in\NN$ and denote by 
$\|\cdot\|_{p,\RD}$ the standard 
norm in $L^p(\RD)$, $1\leq p\leq\infty$.
The Fourier transform $\hat{f}(\xi):=\int_{\RD} e^{-i\xi\cdot{}x}f(x)\,dx$, 
$\xi\in\RD$, of a function $f\in{}L^p(\RD)$ for
$1\leq p\leq2$ is defined in the standard way (see e.g. \cite[Section IV.3]{K}).
 In the case $2<p\leq\infty$
we consider only the functions $f\in{}L^p(\RD)$
whose transforms belong to $L^\pp(\RD)$, 
$\pp:=p/(p-1)$. 
Introduce the difference operator of order $m\in\NN$
acting on functions with domain $\RD$ by \cite[Section~3.3]{Ti}
\begin{equation}
\label{eqfindiff}
     \D_y^m f(x):=\sum_{k=0}^m (-1)^{m-k}
      \Big(\genfrac{}{}{0pt}{}{m}{k}\Big)\,f(x+ky),\qquad x,y\in\RD,
\end{equation}
where $\big(\genfrac{}{}{0pt}{}{m}{k}\big)$ denotes the binomial
coefficient.
Note that the Fourier transform of \eqref{eqfindiff} in $x$ 
equals
$
   \Big(\sum_{k=0}^m (-1)^{m-k}
      \big(\genfrac{}{}{0pt}{}{m}{k}\big)\,e^{iky\cdot\xi}\Big)\,
          \hat{f}(\xi)
  = (e^{iy\cdot\xi}-1)^m\,\hat{f}(\xi),
$
and hence the functions 
\begin{equation}
\label{eqFpair}
        \D_y^mf(x),\qquad (2i\sin(y\cdot \xi/2))^m\,e^{imy\cdot\xi/2}\,\hat{f}(\xi)
\end{equation}
form a Fourier pair for $m\in\NN$. Now we
define $\D_y^mf(x)$ for any $m>0$,
$m\not\in\NN$, as the inverse Fourier transform of the second 
function in \eqref{eqFpair}.
Let $\SD$ denote the unit sphere in $\RD$
and let $dS$ denote the standard 
measure on $\SD$. 
(All the results below with obvious modifications
hold if one replaces
the 
standard measure on $\SD$,
with a 
measure $G$ 
invariant under orthogonal transformations
and supported in the unit ball in $\RD$ 
as in \cite{Cl}.) 
For $1\leq{}p\leq\infty$ set
$$
         \delta_{p,m}[f](y):=\|\Delta_y^m{}f\|_{p,\RD},\qquad
                          y\in\RD,
$$ 
and define a general averaged (integral) \hbox{m.c.} of $f$
as follows. For any $h>0$ in the case $1\leq{}q<\infty$ set 
\begin{equation}
\label{eq200}
   \omega_{p,m,q}[f](h):=\|\d_{p,m}[f](hy)\|_{q,\SD} 
  =  \Big(\int_{\SD} \|\Delta_{hy}^m{}f\|_{p,\RD}^{q}\,dS_y\Big)^{1/q}
 \end{equation}
(where in the case $d=1$ the integral should be interpreted as a sum
over $y\in\SS^0=\{\pm1\}$),
and in the case $q=\infty$ set
$$
        \omega_{p,m,\infty}[f](h):=\sup_{|y|\leq h}\d_{p,m}[f](y)
            =\sup_{|y|\leq h}\|\Delta_y^m{}f\|_{p,\RD}.
$$
The H\"older inequality implies that for any $1\leq q_1\leq q_2\leq\infty$
there exists $c=c(p,q_1,q_2,d)<\infty$ so that 
\begin{equation}
\label{eq0}
   \omega_{p,m,q_1}[f](h) \leq c\,\omega_{p,m,q_2}[f](h),
          \qquad h>0.
\end{equation}
Define the tail integral of the \hbox{F.t.} for $1\leq\pp<\infty$
and for $\pp=\infty$,
respectively, by
\begin{equation}
\label{eq00true}
   \psi_{\pp}[\hat{f}](t)
        := \Big(\int_{|\xi|\geq t} |\hat{f}(\xi)|^{\pp}\,d\xi\Big)^{1/\pp},\qquad
   \psi_{\infty}[\hat{f}](t)
        := \sup_{|\xi|\geq t} |\hat{f}(\xi)|,\qquad t>0.
\end{equation}
Motivated by the results in \cite{Cl}, 
we wish to compare 
the \hbox{m.c.} $\omega_{p,m,q}[f](1/t)$ and the \hbox{F.t.} tail
$\psi_{\pp}[\hat{f}](t)$, as $t\ra\infty$. 
It will be clear from our Theorem \ref{thm1} below
that the natural choice of $q$ for the purpose of such a comparison
is $q=p^\prime$ (see also the discussion preceding Theorem \ref{Besseltail}).
Note next that it is
possible
that  $\psi_{\pp}[\hat{f}](t)$ is rapidly decreasing,
or simply zero, for large $t$ (take e.g., $\hat{f}\in{}C_0^\infty(\RD)$),
whereas the modulus of continuity related
to the $m$th finite difference of $f$ vanishes generally speaking
at the rate $1/t^m$ only.
This motivated the author in \cite{Cl}
to introduce the following modified \hbox{F.t.}
 tails:
For $1\leq\pp<\infty$
\begin{equation}
\label{eq00}
\begin{aligned}
   \psi_{\pp,m}[\hat{f}](t):&= \Big(\int_{\RD} 
   \min\big(1, (|\xi|/t)^{m\pp}\big)
               \,|\hat{f}(\xi)|^{\pp}\,d\xi \Big)^{1/\pp}\\
           &=\Big( t^{-m\pp}
            \int_{|\xi|\leq t} |\xi|^{m\pp}\,|\hat{f}(\xi)|^{\pp}\,d\xi
          +(\psi_{\pp}[\hat{f}](t))^{\pp}\Big)^{1/\pp},\qquad t>0,
\end{aligned}
\end{equation}
and for $\pp=\infty$
$$
          \psi_{\infty,m}[\hat{f}](t):= 
\sup_{\xi\in\RD} \Big(\min\big(1,(|\xi|/t)^{m}\big)\,|\hat{f}(\xi)|\Big),\qquad t>0.
$$ 
It might not be immediately obvious why these tails are useful,
 we give the
reason for that in Theorem~\ref{Besseltail} below.
Note that   
for any $1\leq\pp\leq\infty$, $\psi_{\pp,m}[\hat{f}](t)$
is nonincreasing as $t$ grows,
because so is the function
$\min(1,r/t)$ for any fixed $r>0$. 
Note also that in the case $\hat{f}\in{}C_0^\infty(\RD)$, 
$\hat{f}\not\equiv0$,
there exist $\tilde{c}_1=\tilde{c}_1(\hat{f})>0$, 
$\tilde{c}_2=\tilde{c}_2(\hat{f})>0$ such that
\begin{equation}
\label{eqcinfty}
       \tilde{c}_1\,t^{-m}\leq\psi_{\pp,m}[\hat{f}](t)
            \leq \tilde{c}_2\,t^{-m},\qquad t\ra\infty,
\end{equation}
because $\psi_{\pp}[\hat{f}](t)=0$ for large $t$.

We are ready to state the two results from \cite{Cl} mentioned above.
First \cite[p.~512]{Cl}, for any $d\in\NN$ in the case $2\leq{}p\leq\infty$,
there exists $c_3=c_3(p,d,m)$ such that for for all
functions in the set 
$\{f : f\in{}L^p(\RD)\textrm{ and }\hat{f}\in{}L^{p^\prime}(\RD)\}$,
\begin{equation}
\label{eq1}
    \omega_{p,m,\infty}[f](1/t)\leq c_3\,\psi_{\pp,m}[\hat{f}](t),\qquad t>0.
\end{equation}
Note that in the case $p=2$, \eqref{eq1} holds for all $f\in{}L^2(\RD)$.
(The mentioned formula in \cite{Cl}
involves in fact $\psi_{\pp,m}[\hat{f}](2t)$. The formula \eqref{eq1}
is then true since $\psi_{\pp,m}[\hat{f}](t)$ 
is nonincreasing, and \eqref{eq1} suffices for our purposes.) 
Secondly \cite[(9)]{Cl}, for any $d\in\NN$ in the case $1\leq{}p\leq2$,
for any $a>1$ there exists $c_4=c_4(p,q,d,m,a)$ such that
for any $1\leq{}q\leq\infty$ and for all $f\in{}L^p(\RD)$
\begin{equation}
\label{eq2}
       \psi_{\pp}[\hat{f}](t)\leq c_4\,
          \sum_{k=1}^\infty a^m\omega_{p,m,q}[f](1/(a^kt)),\qquad t>0.
\end{equation}
It would be preferable to have instead of \eqref{eq2}
a formula which does not involve an infinite sum,
i.e., of the type \eqref{eq1}.
It turns out that such a result holds for $d\geq2$ and if $q$ is large enough,
at least $q=\pp$.
\begin{theorem}
\label{thm1}
For any $d\geq2$ and $1\leq{}p\leq2$ there exists $c=c(p,d,m)$ such that
for all $f\in{}L^p(\RD)$
\begin{equation}
\label{eq2p}
    \psi_{\pp,m}[\hat{f}](t) \leq c\,\,\omega_{p,m,\pp}[f](1/t),\qquad t>0.
\end{equation}
\end{theorem}
It is explained in Remark~\ref{remdim1} below
why
\eqref{eq2p} (and even its analog with $\psi_{\pp}[\hat{f}](t)$
in the left-hand side, cf.~\eqref{eqtri} below) fails for $d=1$.
In view of \eqref{eq0}, 
one can replace $\pp$ in the right-hand side in \eqref{eq2p}
with any $q\geq\pp$ (and a different $c$).
It would be interesting to know if one could replace $\pp$
in the right-hand side of \eqref{eq2p} with $1\leq{}q<\pp$.
Note also that for $q\geq\pp$, 
our Theorem \ref{thm1} implies readily
all the statements of \cite[Theorem 2]{Cl},
and is slightly more general since 
no comparison function $s$ as in \cite{Cl} is required.

In the case $p=2$ the following is an immediate consequence 
of \eqref{eq1}, \eqref{eq2p} and \eqref{eq0}.
\begin{corollary}
\label{thm2}
For $d\geq2$ and $p=2$ there exist $c_1,c_2>0$ that depend on $d$
and $m>0$ only, such that for all $f\in{}L^2(\RD)$ and $t>0$
\begin{equation}
\label{eq2p2}
\begin{aligned}
    c_1\,\omega_{2,m,\infty}[f](1/t)\leq \psi_{2,m}[\hat{f}](t)
               &\leq c_2\,\omega_{2,m,2}[f](1/t)\\
   &\leq c_2c(2,2,\infty,d)\,\omega_{2,m,\infty}[f](1/t).
\end{aligned}
\end{equation}
\end{corollary}
The estimates \eqref{eq1}, \eqref{eq2p}, \eqref{eq2p2}
show that that the modified \hbox{F.t.} tail $\psi_{\pp,m}$ 
is more appropriate than the true \hbox{F.t.} tail  $\psi_{\pp}$
to be compared with the \hbox{m.c.} $\omega_{p,m,\pp}$.
The upper estimate in \eqref{eq2p2} is a Jackson-type inequality,
see e.g. \cite[Section I.8]{K}.

It follows from \eqref{eq2p2}
that for $p=2$ an inequality in the direction
opposite to \eqref{eq0} holds: In the case $d\geq2$
for all $f\in{}L^2(\RD)$
\begin{equation}
\label{eq000}
    \omega_{2,m,\infty}[f](1/t)\leq \tilde{c}(2,\infty,2)\,\omega_{2,m,2}[f](1/t),\qquad t>0,
\end{equation}
and so for any $m>0$
all the moduli $\omega_{2,m,q}$, $2\leq{}q\leq\infty$,
are equivalent.
It would be interesting to find a direct proof of \eqref{eq000}.

\begin{figure}[tb]
\begin{center}
\epsfig{file=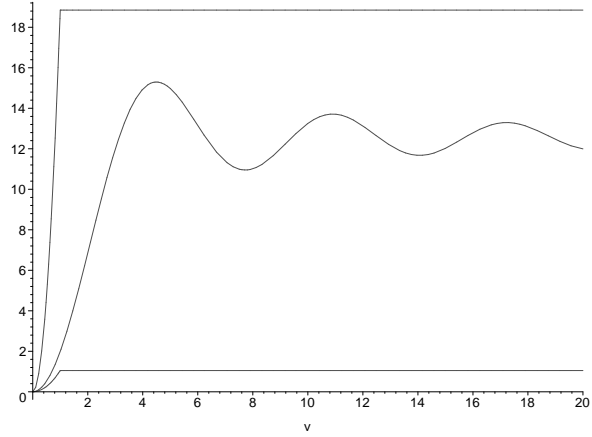,scale=0.4,angle=-90}
\caption{$\frac\pi3\,\min^2(1,v)\leq G_1(v)\leq 6\pi\min^2(1,v)$, 
$0\leq v\leq 20$}
\label{fig}
\end{center}
\end{figure}

We now explain why in the case $d\geq2$ and $1\leq p^\prime<\infty$
the modified \hbox{F.t.}
$\psi_{\pp,m}$ defined in \eqref{eq00}
is a natural quantity to consider. 
For $d\geq2$ and any $0<\alpha<\infty$ define
\begin{equation}
\label{eq_Besselt}
   G_\alpha(|w|):=
         2^\alpha\int_{\SD}\big(1-\cos (y\cdot{}w)\big)^\alpha\,dS_y,
      \qquad w\in\RD.
\end{equation}
Recalling that the functions in \eqref{eqFpair} form a Fourier pair
and using the Hausdorff--Young 
inequality (see the proof of Theorem \ref{thm1}
below and the proof of \eqref{eq1} in \cite{Cl})
we can compare for $y\in\SD$ fixed and $t>0$,
the $L^p$ norm of $\Delta_y^m{}f$, $\delta_{p,m}[f](y)$,
and the $L^{p^\prime}$ norm of the function 
$(2i\sin(y\cdot \xi/(2t)))^m\,e^{imy\cdot\xi/2}\,\hat{f}(\xi)$
(in the case $p>2$ we assume in addition
as before that $\hat{f}\in L^{p^\prime}(\RD))$. 
Raising both quantities to the power $q=p^\prime$
(this explains why the choice $q=p^\prime$ is natural)
we can compare the quantities
$(\omega_{p,m,p^\prime}[f](t))^{\pp}$
and 
$$
   \int_{\RD}
        2^{m\pp/2}
    \bigg(\int_{\SD}\Big(1-\cos\frac{y\cdot\xi}{t}\Big)^{m\pp/2}\,dS_y\bigg)
                      |\hat{f}(\xi)|^{\pp}\,d\xi.
$$
With this in mind, for $d\geq2$ and
any $\hat{f}\in L^{\pp}(\RD)$, $1\leq p<\infty$,
 we introduce the {\em Bessel tail\ }of the Fourier transform
$$
      \Psi_{\pp,m}[\hat{f}](t)
          :=
   \bigg(\int_{\RD} G_{m\pp/2}\big({|\xi|}/{t}\big) 
                      |\hat{f}(\xi)|^{\pp}\,d\xi\bigg)^{1/\pp}
$$
and $\Psi_{\infty,m}[\hat{f}](t):=\psi_{\infty,m}[\hat{f}](t)$.
From the above discussion, $\omega_{p,m,p^\prime}[f](t)$
can be compared with the Bessel tail $\Psi_{\pp,m}[\hat{f}](t)$
as in \eqref{eq1}, \eqref{eq2p}.
The relevance of the modified \hbox{F.t.}
$\psi_{\pp,m}$ is now apparent from the following
\begin{theorem}\label{Besseltail}
For any $d\geq2$ and any $0<\alpha<\infty$, 
there exist $C_1,C_2>0$ that depend only on $d$ and $\alpha$
so that 
\begin{equation}
\label{eqBessel}
      C_1(\min(1,v))^{2\alpha}\leq G_\alpha(v)
          \leq C_2(\min(1,v))^{2\alpha},\qquad v\geq0,
\end{equation}
and hence for some $\tilde{C}_1,\tilde{C}_2>0$ that depend
on $1\leq p<\infty$, $m>0$, $d$ only,
$$
      \tilde{C}_1\Psi_{\pp,m}[\hat{f}](t)
                      \leq \psi_{\pp,m}[\hat{f}](t)
            \leq\tilde{C}_2\Psi_{\pp,m}[\hat{f}](t)
$$
for all $\hat{f}\in L^{\pp}(\RD)$ and all $t>0$.
\end{theorem}
The relation \eqref{eqBessel} is illustrated for 
$\alpha=1$, $C_1=\pi/3$, $C_2=6\pi$ in Fig.~\ref{fig}
(note that $G_1(v)=4\pi\big(1-\frac{\sin v}{v})$).

We now describe the relationship between
the true and the modified \hbox{F.t.} tails,
 $\psi_\pp$ and  $\psi_{\pp,m}$, respectively.
From \eqref{eq00} it is clear that
for any $1\leq\pp\leq\infty$, $m>0$, and $\hat{f}\in{}L^\pp(\RD)$
\begin{equation}
\label{eqtri}
         \psi_\pp[\hat{f}](t)\leq\psi_{\pp,m}[\hat{f}](t),\qquad t>0.
\end{equation}
The following statement gives a 
converse to \eqref{eqtri} that is optimal on the power scale
(see Remark~\ref{remstar} below).
\begin{theorem}
\label{lem4.1}
Let $d\in\NN$, $m>0$, $\alpha>0$, $1\leq\pp\leq\infty$, 
$g{}\in{}L^\pp(\RD)$. All the constants below depend on $g$
and are strictly positive and finite.

1. Let $1\leq\pp<\infty$. 

(i) If $\psi_\pp[g{}](t)\leq c_2(g)\cdot{}t^{-\alpha}$, $t\ra\infty$, 
then as $t\ra\infty$,  
$$
\psi_{\pp,m}[g{}](t)\leq b_2(g)
      \cdot\begin{cases}
                                      t^{-\alpha},&0<\alpha<m\cr
                                   t^{-\alpha}\,(\log{}t)^{1/\pp},&\alpha=m\cr
                             t^{-m},&\alpha>m.
      \end{cases}
$$

(ii) If $c_1(g)\cdot{}t^{-\alpha}\leq
\psi_{\pp,m}[g{}](t)\leq c_2(g)\cdot{}t^{-\alpha}$, $t\ra\infty$, 
then as $t\ra\infty$, 
$$
\psi_{\pp}[g{}](t)\geq b_1(g)
      \cdot\begin{cases}
                                      t^{-\alpha},&0<\alpha<m\cr
                                      0,&\alpha\geq m.
      \end{cases}
$$

2. Let $\pp=\infty$.

(i) If $\psi_\infty[g{}](t)\leq c_2(g)\cdot{}t^{-\alpha}$, $t\ra\infty$, 
then as $t\ra\infty$,  
$$
\psi_{\infty,m}[g{}](t)\leq b_2(g)
      \cdot\begin{cases}
                                    t^{-\alpha},&0<\alpha< m\cr
                                    t^{-m},&\alpha\geq m.
      \end{cases}
$$

(ii) If $c_1(g)\cdot{}t^{-\alpha}\leq
\psi_{\infty,m}[g{}](t)\leq c_2(g)\cdot{}t^{-\alpha}$, $t\ra\infty$, 
then as $t\ra\infty$, 
$$
\psi_{\infty}[g{}](t)\geq b_1(g)
      \cdot\begin{cases}
                                      t^{-\alpha},&0<\alpha<m\cr
                                      0,&\alpha\geq m.
      \end{cases}
$$
\end{theorem}
%

We state finally a result that was applied in a study of the scaling
of entanglement entropy for a certain physical system in \cite{GiKl,Gi2}.
\begin{corollary}
\label{CH3_s3.lem.1} 
Let $d\in\NN$.
Assume that ${}f\in{}L^2(\RD)$. 
Then for some $c_1=c_1(f)>0$, $c_2=c_2(f)>0$
and some $\gamma=\gamma(f)\in(0,2)$, $f$ satisfies
\begin{equation}
\label{eqneww}
        c_1\,\eps^\gamma\leq
   \int_{\SD}\|{}f(\cdot+\eps{}y)-{}f(\cdot)\|_{2,\RD}^2 
\,dS_y\leq c_2\,\eps^\gamma,
             \qquad0\leq\eps\leq1,
\end{equation}
if and only if there exist $b_1=b_1(f)>0$, $b_2=b_2(f)>0$ such that
\begin{equation}
\label{eqnewww}
          b_1\,t^{-\gamma} \leq 
                   \int_{|\xi|\geq{}t} 
                      |\hat{f}(\xi)|^2\,d\xi 
           \leq b_2\,t^{-\gamma},\qquad t\geq1.
\end{equation}
\end{corollary}
Note that Corollary \ref{CH3_s3.lem.1} is true in all dimensions:
In the proof of Corollary \ref{CH3_s3.lem.1} below, we consider
the cases $d\geq2$ and $d=1$ separately. In the former case we employ the 
general results stated above. In the case $d=1$ we give a direct proof
using in particular the ideas in the proofs 
of \cite[Lemma 2.10, Lemma 4.2]{BCT}. The main reason why,
in the case $d=1$,
Corollary \ref{CH3_s3.lem.1} is true despite the fact
that Theorem \ref{thm1} fails, is because of the explicit (power-type)
form
of the estimates in Corollary \ref{CH3_s3.lem.1}.

The equivalence of the {\em upper\ }estimates in \eqref{eqneww}
and in \eqref{eqnewww} is well-known: it follows
follows e.g. from \cite[Lemma~3.3.1]{MS},
and also from the results
obtained in \cite{Cl}. Note that \cite[Lemma~3.3.1]{MS} deals 
with a Besov space $B_{2,\infty}^{s}(\RD)$, $s>0$, $s\not\in\NN$
(the case $s=\gamma/2\in(0,1)$ is relevant for the upper estimates 
in Corollary~\ref{CH3_s3.lem.1}).
In \cite{BCT}, \eqref{eqnewww} is derived in the case $\gamma=1$
from a more restrictive pointwise condition
\begin{equation}
\label{eq_ts}
      c_1\,|y|^\gamma\leq(\d_{2,1}[{}f](y))^2
               \leq c_2\,|y|^\gamma,\qquad|y|\leq1,
\end{equation}
which in general can not be reversed because $\d_{2,1}[{}f](y)$ 
can have singular directions.

A simple example of a function that 
satisfies \eqref{eqneww} with $\gamma=1$
is the characteristic function of a compact
set with $C^1$ boundary.
For any $0<\gamma<1$ there exsits a compact set whose 
characteristic function satisfies \eqref{eqneww}, see \cite[Lemma 2.9]{Gi2}.

It turns out that Corollary~\ref{CH3_s3.lem.1} fails for $\gamma=2$,
see Remark~\ref{remthf}.
Note that for $\gamma>2$ the condition \eqref{eqneww} is not satisfied
for any $f\in{}C_0^\infty(\RD)$, $f\not\equiv0$,
 because it involves the finite difference
of order~$1$.

Theorems \ref{thm1}, \ref{Besseltail}, \ref{lem4.1}
and Corollary \ref{CH3_s3.lem.1} 
are proved in Section \ref{sec1.3}.
%
%
%
%
\section{Proofs and concluding remarks}
\label{sec1.3}
\begin{proof}[Proof of Theorem~\ref{thm1}]
First consider the case $p=1$, $\pp=\infty$
(in this case the proof below goes through for all $d\in\NN$). 
We have to prove that
for some $c=c(p,d,m)$,
$\psi_{\infty,m}[\hat{f}](t)\leq c\,\omega_{1,m,\infty}[f](1/t)$, $t>0$. 
Recall that $\|\hat{f}\|_{\infty,\RD}\leq\|f\|_{1,\RD}$.
Since the functions \eqref{eqFpair}
form a Fourier pair we have
\begin{equation}
\label{eq30}
\begin{aligned}
     \omega_{1,m,\infty}[f](1/t) &= \sup_{|y|\leq1/t}\d_{1,m}(y) 
             \geq \sup_{|y|\leq1/t} \sup_{\xi\in\RD}
                \big(|2i\sin (y\cdot \xi/2)\big|^m\,|\hat{f}(\xi)|\big)\\
    &= 2^m\,\sup_{\xi\in\RD}\Big( |\hat{f}(\xi)| 
            \sup_{|z|\leq1} \big|\sin (z\cdot \xi/(2t))\big|^m \Big).
\end{aligned}
\end{equation}
It is an easy exercise to prove
that for any $\xi\in\RD$
$$
\begin{aligned}
      \sup_{|z|\leq1} \big|\sin (z\cdot \xi/(2t))\big| &=
      \begin{cases} 1, &|\xi|\geq |t\pi|\cr
                            \sin(|\xi|/(2t)), &|\xi|< |t\pi|
     \end{cases}\\
   &\geq\frac14\,\min(1,|\xi|/t).
\end{aligned}
$$
This together with \eqref{eq30} and \eqref{eq00}
proves the result.

Consider now the case $d\geq2$, $1<p\leq2$, $\pp<\infty$.
By the Hausdorff--Young inequality
$\|\hat{f}\|_{\pp,\RD}\leq(2\pi)^{d/\pp}\|f\|_{p,\RD}$,
$1<p\leq2$, and using the fact that \eqref{eqFpair} 
is a Fourier pair we get as in \cite[(4.7)]{Cl} 
\begin{equation}
\label{eq300}
\begin{aligned}
   (\d_{p,m}[f](y/t))^{\pp}
    &\geq (2\pi)^{-d}\int_{\RD} \big|2i\sin (y\cdot \xi/(2t))\big|^{m\pp}\,
             |\hat{f}(\xi)|^{\pp}\,d\xi\\
   &= (2\pi)^{-d}\,2^{m\pp/2}\int_{\RD} |\hat{f}(\xi)|^{\pp}
        \,\Big(1-\cos\frac{y\cdot\xi}{t}\Big)^{m\pp/2}\,d\xi.
\end{aligned}
\end{equation}
Integrating over $y$ and recalling \eqref{eq_Besselt} we obtain
\begin{equation}
\label{eq3}
\begin{aligned}
    (\omega_{p,m,\pp}[f](1/t))^{\pp} 
       &=\int_{\SD}(\d_{p,m}[f](y/t))^{\pp}\,dS_y\\
   &\geq (2\pi)^{-d}\int_{\RD}G_{m\pp/2}(|\xi|/t)\, |\hat{f}(\xi)|^{\pp}\,d\xi.
\end{aligned}
\end{equation}
Now we need the following elementary result, the proof is given
after the end of the present proof.
\begin{lemma}
\label{lem1}
Let $d\geq2$ and fix any $0<\alpha<\infty$.
The function $G_\alpha(v)$, $v\geq0$,
defined in \eqref{eq_Besselt}
satisfies the following:
For any $v_0>0$
there exist $C_\alpha(v_0),c_\alpha(v_0)>0$ such that
\begin{equation}\label{eqlem101}
       {}G_\alpha(v)\geq C_\alpha(v_0) >0,\qquad v\geq v_0
\end{equation}
and
\begin{equation}\label{eqlem102}
       G_\alpha(v)\geq c_\alpha(v_0)\, v^{2\a}, \qquad 0\leq v\leq v_0.
\end{equation}
\end{lemma}
Now we write the integral in \eqref{eq3} as a sum of 
$\int_{|\xi|\leq{}t}$ and $\int_{|\xi|\geq{}t}$, apply
Lemma~\ref{lem1} with $\alpha=m\pp/2$, $v_0=1$, and 
set $c_1(p,d,m):=c(m,d)\cdot\min\big( c_\alpha(1),C_\alpha(1))$
to obtain 
\begin{equation}
\label{eq5}
\begin{aligned}
    (\omega_{p,m,\pp}[f](1/t))^{\pp} 
       &\geq c_1(p,d,m)\,
        \Big( \int_{|\xi|\leq{}t}(|\xi|/t)^{m\pp}\,|\hat{f}(\xi)|^{\pp}\,d\xi 
          + \int_{|\xi|\geq{}t}|\hat{f}(\xi)|^{\pp}\,d\xi\Big)\\
      &=c_1(p,d,m)\,
             (\psi_{\pp,m}[\hat{f}](t))^{\pp}
\end{aligned}
\end{equation}
which finishes the proof of Theorem \ref{thm1}.
\end{proof}
\begin{remark}\label{remdim1}
We show now that \eqref{eq2p} is not true in the case $d=1$ and $p=2$.
More precisely, we show that for $d=1$ 
there is no $c\in(0,+\infty)$ so that for all $f\in{}L^2(\RR)$
\begin{equation}\label{eqnn}
    \psi_{2}[\hat{f}](t) \leq c\,\,\omega_{2,1,2}[f](1/t),\qquad t>0.
\end{equation}
By \eqref{eqFpair} with $m=1$ using the Parseval formula
we find
\begin{equation}\label{eqParseval}
          (\omega_{2,1,2}[f](1/t))^2 = \frac4\pi\int_{\RR}
                   \sin^2\bigg(\frac\xi{2t}\bigg)\,|\hat{f}(\xi)|^2\,d\xi.
\end{equation}
For $\xi\in\RR$, $t,\tilde{c}>0$, introduce the notation
$$
      H(\xi,t,\tilde{c}) = \sin^2\bigg(\frac\xi{2t}\bigg) - \tilde{c}
              \,\sgn\bigg[\Big(\frac{|\xi|}t-1\Big)_+\bigg]
$$
where $\sgn a=1,0,-1$ for $a>0$, $a=0$, $a<0$, respectively,
and $(a)_+=\max(0,a)$.
If \eqref{eqnn} were true for some $c\in(0,+\infty)$ then 
there would exist $\tilde{c}\in(0,+\infty)$ such that
$$
         \int_{\RR} H(\xi,t,\tilde{c}) \,|\hat{f}(\xi)|^2\,d\xi \geq0\qquad
       \textrm{for all } f\in L^2(\RR),\quad t>0.
$$
But 
$$
\begin{aligned}
    \{\, |\hat{f}|^2 \,:\,f\in L^2(\RR)\, \} &= 
          \{ \,|f|^2 \,:\,f\in L^2(\RR) \,\}\\
      &=
         \{\, g \,:\,g\in L^1(\RR)\textrm{ and }g\geq0
          \textrm{ almost everywhere (a.e.)} \, \}.
\end{aligned}
$$
We have arrived at a contradiction: If
$$
       \int_{\RR} H(\xi,t,\tilde{c}) \,g(\xi)\,d\xi \geq0\qquad
       \textrm{for all } g\in L^1(\RR),\,g\geq0\textrm{ a.e., } t>0,
$$
then we must have $H(\xi,t,\tilde{c})\geq0$ for a.e.~$\xi$
and all $t>0$, which is clearly false for any choice of $\tilde{c}>0$.
This proves the result. 

We note that the basic reason
for the inapplicability of Theorem \ref{thm1} to the case $d=1$
is that for any $v_0>0$ it is not possible
to insert a constant function between
the graph of $ \sin^2 v$ and the real axis on
the interval $[v_0,+\infty)$, cf.~Fig.~\ref{fig}.
(In the case $d=1$, $\SS^0=\{\pm1\}$ and so
the function $G_\alpha$ in \eqref{eq_Besselt}
would be given by $2^{2\alpha+1}\big(\sin^2(w/2)\big)^\alpha$,
$w\in\RR$, and there is no helpful averaging over $\SS^{d-1}$.)
\end{remark}
\begin{proof}[Proof of Lemma \ref{lem1}]
Recall that $d\geq2$.
If $d\geq3$ then integrating over the $d-2$ angles in \eqref{eq_Besselt}
as in \cite[(II.3.4.2)]{GS} 
we obtain for $v\geq0$
\begin{equation}
\label{eq101z}
    G_\alpha(v) = |\SS^{d-2}|\int_{0}^\pi (1-\cos(v\cos\theta))^\a
 \,\sin^{d-2}\theta\,d\theta.
\end{equation}
If $d=2$ then \eqref{eq101z} and the subsequent formlulae
hold true with the convention $|\SS^{0}|=2$.

We show first there is $M_\alpha>0$
such that ${}G_\alpha(v)\geq{}M_\alpha>0$ for sufficiently large $v$. 
Indeed, in the case $\alpha\geq1$ by the H\"older inequality
there is $\tilde{C}_\alpha>0$ so that
\begin{equation}
\label{eq101}
\begin{aligned}
   \int_{0}^\pi \big( (1-\cos(v\cos\theta))
 \,&\sin^{(d-2)/\a}\theta\big)^\a\,d\theta \\
         &\geq \tilde{C}_\alpha\,\Big(\int_{0}^\pi (1-\cos(v\cos\theta))
 \,\sin^{(d-2)/\a}\theta\,d\theta \Big)^\a \\
   &\geq \tilde{C}_\alpha\,\Big(\int_{0}^\pi (1-\cos(v\cos\theta))
      \,\sin^{d-2}\theta\,d\theta \Big)^\a 
\end{aligned}
\end{equation}
where the second inequality follows from
\begin{equation}
\label{eqeasy}
       w^\beta \geq w,\qquad 0<\beta\leq1,\qquad  0\leq w\leq1
\end{equation}
with $w=\sin^{d-2}\theta$ and $\beta=1/\a$.
In the case $0<\a<1$, in view of \eqref{eqeasy}
with $w=(1-\cos(v\cos\theta))/2$ and $\beta=\a$ we obtain
\begin{equation}
\label{eq101p}
\begin{aligned}
   \int_{0}^\pi  (1-\cos(v\cos\theta))^\a
 \,\sin^{d-2}\theta\,d\theta 
          &= 2^{\a}\,\int_{0}^\pi 
       \Big(\frac{1-\cos(v\cos\theta)}{2}\Big)^\a
            \,\sin^{d-2}\theta\,d\theta \\
        &\geq 2^{\a-1}\,\int_{0}^\pi  (1-\cos(v\cos\theta))
            \,\sin^{d-2}\theta\,d\theta.
\end{aligned}
\end{equation}
Note that
\begin{equation}
\label{eq1019}
\begin{aligned}
     |\SS^{d-2}|\int_{0}^\pi (1-\cos(v\cos\theta))
      \,\sin^{d-2}\theta\,d\theta
  &= \int_{\SD}\Big(1-\cos\big((v,0,\cdots,0)\cdot y\big)\Big)\,dS_y\\
   &=  |\SD|\,\big(1 - 2^s\,\Gamma(s+1)\,
          v^{-s}\,J_s(v)\big)
\end{aligned}
\end{equation}
where $s=(d-2)/2$, $J_s$ is the Bessel function,
and we have used \cite[(II.3.4.2)]{GS} and \cite[(8.411.4)]{GrRy}.
By \cite[(8.451)]{GrRy}, $v^{-s}J_s(v)=O(v^{-s-1/2})\ra0$,
as $v\ra\infty$ for $d\geq2$. Hence the right-hand side of \eqref{eq1019},
and also of \eqref{eq101} and \eqref{eq101p}, tends to a strictly positive limit,
as $v\ra\infty$.
Therefore for any fixed $0<\alpha<\infty$, there exist $M_\alpha>0$
and $v_1(\alpha)$ so that
${}G_\alpha(v)\geq{}M_\alpha>0$ for $v\geq{}v_1(\alpha)$. 

But ${}G_\alpha(v)$ does not have  zeros other than $v=0$.
Since $G_\alpha$ is continuous, it is for any $v_0>0$
bounded away from zero on the compact $[v_0,v_1(\alpha)]$. 
This proves \eqref{eqlem101}.

Let us now prove \eqref{eqlem102} for $v_0>0$ small enough.
We can rewrite \eqref{eq101z} as 
 \begin{equation}
\label{eqlem103}
    G_\alpha(v) = 2^{\alpha}\,|\SS^{d-2}|\int_{0}^\pi 
         \sin^{2\alpha}(v\cos\theta)
 \,\sin^{d-2}\theta\,d\theta.
\end{equation}
Using the elementary estimate
\begin{equation}\label{eqlem104}
         \sin x \geq \frac2\pi\, x,\qquad 0\leq x\leq\frac\pi2,
\end{equation}
we conclude that 
$$
    G_\alpha(v) \geq v^{2\alpha}\cdot c(\alpha,d) \int_{0}^\pi 
         \cos^{2\alpha}\theta
 \,\sin^{d-2}\theta\,d\theta, \qquad 0\leq v\leq \frac\pi2,
$$
which proves \eqref{eqlem102} for any $0<v_0\leq\pi/2$.
But now if we take any $v_0>\pi/2$ then
using \eqref{eqlem102} for $0\leq v\leq\pi/2$ and 
the fact that $G_\alpha(v)$ is bounded away from zero
for $\pi/2\leq v\leq v_0$ by \eqref{eqlem101} we can 
always find $c_\alpha(v_0)>0$ small enough so that \eqref{eqlem102}
holds for $0\leq v\leq v_0$.
The proof of Lemma~\ref{lem1} is complete.
\end{proof}
\begin{proof}[Proof of Theorem \ref{Besseltail}]
Recall that $d\geq2$.
The lower estimate in \eqref{eqBessel} follows immediately
 from Lemma \ref{lem1} with $v_0=1$.
As for the upper estimate, we note first that by the definition
\eqref{eq_Besselt}, the function $G_\alpha(v)$ is bounded 
above for $v\geq1$.
Next, using the estimate $\sin x\leq x$, $0\leq x\leq1$,
in place of \eqref{eqlem104}, we derive from \eqref{eqlem103}
the upper estimate
$$
           G_\alpha(v) \leq v^{2\alpha}\cdot \tilde{c}(\alpha,d) \int_{0}^\pi 
         \cos^{2\alpha}\theta
                    \,\sin^{d-2}\theta\,d\theta, \qquad 0\leq v\leq 1.
$$
This proves the upper estimate in \eqref{eqBessel}.
The proof of Theorem \ref{Besseltail} is complete.
\end{proof}
\begin{proof}[Proof of Theorem \ref{lem4.1}]
Recall that $d\in\NN$,  $0<\alpha<\infty$, $m>0$.

1. Consider first the case $1\leq\pp<\infty$.

(i) Assume $\psi_{\pp}[{}g{}](t)\leq{}c_1(g)\cdot{}t^{-\alpha}$, $t\geq1$.
By \eqref{eq00}
\begin{equation}
\label{eq201}
 (\psi_{\pp,m}[{}g{}](t))^\pp=   t^{-m\pp}
            \int_{|\xi|\leq t} |\xi|^{m\pp}\,|{}g{}(\xi)|^\pp\,d\xi
          +(\psi_{\pp}[{}g{}](t))^\pp.
\end{equation}
Then the following gives the result for the case 1(i), all values of $\alpha$:
$$
\begin{aligned}
   t^{-m\pp} &\int_{|\xi|\leq t} |\xi|^{m\pp}\,|{}g{}(\xi)|^\pp\,d\xi\\
       &\leq t^{-m\pp}\Big(\int_{|\xi|\leq1}|{}g{}(\xi)|^\pp\,d\xi
                    + \sum_{k=0}^{[\log_2t]} 
             \int_{2^{k}\leq|\xi|\leq2^{k+1}}|\xi|^{m\pp}\,
                    |{}g{}(\xi)|^\pp\,d\xi\Big)\\
          &\leq t^{-m\pp}\Big(\const(g)+ c_1(g)\sum_{k=0}^{[\log_2t]} 
             \big(2^{(m-\alpha)\pp}\big)^k\Big)
\end{aligned}
$$
where $[\cdot]$ denotes the integer part of a real number.
It is explained in Remark \ref{remstar} below
why the order in $t$ cannot be improved
 in the case 1(i) and in all other cases.

(ii) Assume $c_1(g)\cdot t^{-\alpha}\leq
\psi_{\pp,m}[{}g{}](t)\leq c_2(g)\cdot t^{-\alpha}$, 
$t\geq1$ ($\pp<\infty$).
The example ${}g{}\in{}C_0^\infty(\RD)$, ${}g{}\not\equiv0$,
shows that in the case $\alpha\geq{}m$ we can only claim
the trivial bound $\psi_\pp[{}g{}](t)\geq0$, $t\ra\infty$.
Let now $0<\alpha<m$. We use the idea in the proof of 
\cite[Lemma~4.2]{BCT}. 
Let $0<B\leq1$ be a number to be chosen later.
By the definition \eqref{eq201}
\begin{equation}\label{eqlem120}
    (\psi_{\pp,m}[{}g{}](t))^\pp = 
           \Big( \int_{|\xi|\leq Bt} + \int_{|\xi|\geq Bt} 
       \Big) \min^{m\pp}(1,|\xi|/t)\,|{}g{}(\xi)|^\pp\,d\xi.
\end{equation}
We have
$$
\begin{aligned}
   \int_{|\xi|\leq Bt} &\min^{m\pp}(1,|\xi|/t)\,|{}g{}(\xi)|^\pp\,d\xi
  =  t^{-m\pp} \int_{|\xi|\leq Bt} |\xi|^{m\pp}\,|{}g{}(\xi)|^\pp\,d\xi\\
       & \leq t^{-m\pp} \bigg(\int_{|\xi|\leq B}|\xi|^{m\pp}\,|{}g{}(\xi)|^\pp\,d\xi
   + \sum_{l=1}^{[\log_2 t]+1} 
             \int_{2^{l-1}B\leq|\xi|\leq2^{l}B}
                            |\xi|^{m\pp}\,|{}g{}(\xi)|^\pp\,d\xi
    \bigg)\\
  & \leq t^{-m\pp} \bigg(\|g\|_{\pp,\RD}^{\pp}
     + \sum_{l=1}^{[\log_2 t]+1} 
            (2^{l}B)^{m\pp}\cdot  c_2^{\pp}(g)\cdot (2^{l-1}B)^{-\alpha\pp}
    \bigg)\\
       &\leq t^{-m\pp}\cdot\const(g)\cdot t^{(m-\alpha)\pp}
             \cdot B^{(m-\alpha)\pp},\qquad t\to\infty.
\end{aligned}
$$
Therefore choosing $0<B\leq1$ small enough (recall that $m-\alpha>0$)
we obtain that
\begin{equation}\label{eqlem121}
        \int_{|\xi|\leq Bt}\min^{m\pp}(1,|\xi|/t)\,|{}g{}(\xi)|^\pp\,d\xi
             \leq \frac{c_1^{\pp}(g)}2 \cdot t^{-\alpha\pp},\qquad t\to\infty.
\end{equation}
Since $\psi_{m,\pp}[g](t)\geq c_1^{\pp}(g)\cdot t^{-\alpha\pp}$,
$t\to\infty$,
\eqref{eqlem121} together with  \eqref{eqlem120} gives
$$
   \int_{|\xi|\geq Bt} 
     \min^{m\pp}(1,|\xi|/t)\,|{}g{}(\xi)|^\pp\,d\xi
      \geq \frac{c_1^{\pp}(g)}2 \cdot t^{-\alpha\pp},\qquad t\to\infty. 
$$
Setting $s=Bt$ we obtain
\begin{equation}\label{eqlem110prim}
    \int_{|\xi|\geq s} 
     \min^{m\pp}(1,B|\xi|/s)\,|{}g{}(\xi)|^\pp\,d\xi
      \geq {c_1(g)}B^{\alpha\pp}\cdot s^{-\alpha\pp},
              \qquad s\to\infty.
\end{equation}
Noting that 
$
    1\geq\min(1,B|\xi|/s)
$
we conclude from \eqref{eqlem110prim} that
$$
   \psi_{\pp}^\pp[g](t) = \int_{|\xi|\geq s} 
     |{}g{}(\xi)|^\pp\,d\xi
               \geq {c_1(g)}B^{\alpha\pp}\cdot s^{-\alpha},\qquad s\to\infty.
$$

2. The case $\pp=\infty$.

(i) Assume $\psi_\infty[{}g{}](t)\leq{}c_1(g)\cdot t^{-\alpha}$, $t\geq1$.
Note that 
\begin{equation}
\label{eqggg}
        \psi_{\infty,m}[{}g{}](t) 
              = \max\Big(t^{-m}\sup_{|\xi|\leq t} |\xi|^m|{}g{}(\xi)|,
                          \,\sup_{|\xi|\geq t} |{}g{}(\xi)|\Big).
\end{equation}
Clearly
\begin{equation}
\label{eqgg}
       |{}g{}(\xi)|\leq \sup_{|\eta|\geq|\xi|} |{}g{}(\eta)|
             \leq c_1(g)\cdot |\xi|^{-\alpha},\qquad|\xi|\geq1.
\end{equation}
Therefore
$$
\begin{aligned}
    \sup_{|\xi|\leq t} |\xi|^m&|{}g{}(\xi)|
        =\max\big(
             \sup_{|\xi|\leq 1} |\xi|^m|{}g{}(\xi)|,
                \sup_{1\leq|\xi|\leq t} |\xi|^m|{}g{}(\xi)|\big)\\
       &\leq \max\Big(c({}g{}),\,
               c_1({}g{})\cdot \sup_{1\leq|\xi|\leq{}t}|\xi|^{m-\alpha}\Big)
          =\begin{cases}C({}g{}),&0<\alpha<m\cr
                        C_1({}g{})\cdot{}t^{m-\a},&\a\geq m.
             \end{cases}
\end{aligned}
$$
This together with \eqref{eqggg} proves the case 2(i).

(ii) Let $c_1(g)\cdot t^{-\alpha}\leq 
\psi_{\infty,m}[{}g{}](t)\leq c_2(g)\cdot t^{-\alpha}$, $t\to\infty$.
Again if $\alpha\geq{}m$ then the example of ${}g{}\in{}C_0^\infty(\RD)$,
$g\not\equiv0$, shows that generally speaking only the trivial
bound $\psi_{\infty}[{}g{}](t)\geq0$ holds for large $t$. 
Let now $0<\alpha<m$.
Let $0<B\leq1$ be a number to be chosen later.
We have
$$
\begin{aligned}
   \psi_{\infty,m}[{}g{}](t) 
              &= \sup_{\xi\in\RD}
         \min^m(1,|\xi|/t)\,|g(\xi)|\\
   &=\max\Big( \sup_{|\xi|\leq Bt}
         \min^m(1,|\xi|/t)\,|g(\xi)|,\,\sup_{|\xi|\geq Bt}
         \min^m(1,|\xi|/t)\,|g(\xi)|
         \Big).
\end{aligned}
$$
Now 
$$
\begin{aligned}
      \sup_{|\xi|\leq Bt}
         &\min^m(1,|\xi|/t)\,|g(\xi)| 
   =    
              t^{-m} \sup_{|\xi|\leq Bt} |\xi|^m\,|g(\xi)|\\
    &\leq    t^{-m} \max\bigg( \sup_{0\leq|\xi|\leq B} |\xi|^m\,|g(\xi)|,
         \max_{l=1,\cdots,[\log_2t]+1}
          \sup_{2^{l-1}B\leq|\xi|\leq 2^lB} |\xi|^m\,|g(\xi)|\bigg)\\
    &\leq    t^{-m} \max\bigg( B^m\|g\|_{\infty,\RD},
         \max_{l=1,\cdots,[\log_2t]+1}
          2^{lm} B^m\cdot c_2(g)\cdot (2^{l-1}B)^{-\alpha}\bigg)\\
&\leq \const(g)\cdot B^{m-\alpha}\cdot t^{-\alpha},\qquad t\to\infty.
\end{aligned}
$$
Choosing $B$ small enough (note that $m-\alpha>0$) we obtain
$$
  {c_1(g)}\cdot t^{-\alpha} \leq
     \psi_{\infty,m}[{}g{}](t) 
         \leq \max\Big( \frac{c_1(g)}2\cdot t^{-\alpha},\,\sup_{|\xi|\geq Bt}
         \min^m(1,|\xi|/t)\,|g(\xi)|
         \Big),\qquad t\to \infty,
$$
which implies that
$\sup_{|\xi|\geq Bt}
         \min^m(1,|\xi|/t)\,|g(\xi)| \geq {c_1(g)}\cdot t^{-\alpha}$.
Replacing $s=Bt$ we obtain
\begin{equation}\label{eqlem110}
     \sup_{|\xi|\geq s}
         \min^m(1,B|\xi|/s)\,|g(\xi)| \geq {c_1(g)}B^\alpha\cdot s^{-\alpha},
              \qquad s\to\infty.
\end{equation}
Noting that 
$
    1\geq\min(1,B|\xi|/s)
$
we conclude from \eqref{eqlem110} that
$$
   \psi_\infty[g](t) = \sup_{|\xi|\geq s} |g(\xi)|
      \geq {c_1(g)}B^\alpha\cdot s^{-\alpha},\qquad s\to\infty.
$$
This completes the proof of Theorem \ref{lem4.1}.
\end{proof}
\begin{remark}
\label{remstar}
Let us explain why the estimates in Theorem~\ref{lem4.1}
have the best possible order in $t$. 
In view of \eqref{eqtri}, the estimates 1(i) and 2(i) 
for $0<\a<m$ can not be improved.
The example of $g\in{}C_0^\infty(\RD)$, $g\not\equiv0$,
in view of the lower estimate in \eqref{eqcinfty} 
shows that estimates 1(i) for $\a>m$ 
and 1(ii), 2(i), 2(ii) for $\a\geq{}m$ can not be improved.
The example of $g\in{}L^\pp(\RD)$, $1\leq\pp<\infty$, 
such that $g(\xi)=|\xi|^{-(d/\pp)-\a}$ for $|\xi|\geq1$ and 
$g$ smooth for $|\xi|<1$ shows that the estimate 1(ii)
for $0<\a<m$ can not be improved 
(note that the mentioned $g$
satisfies $\psi_{\pp,m}[g](t)]\leq
c_1(g)\cdot{}t^{-\a}$, $t\ra\infty$,
and also $\psi_\pp[g](t)=
\big(|\SD|/(\a\pp)\big)^{1/\pp}\cdot{}t^{-\a}$, $t\geq1$).
Choosing $g\in{}L^\infty(\RD)$
such that $g(\xi)=|\xi|^{-\a}$ for $|\xi|\geq1$ and 
$g$ smooth for $|\xi|<1$ shows that the estimate 2(ii)
for $0<\a<m$ can not be improved 
(for this $g$
we have $\psi_\infty[g](t)=t^{-\a}$, $t\geq1$).
Finally, the example ${}g{}\in{}L^\pp(\RD)$, $1\leq\pp<\infty$,
such that ${}g{}(\xi)=|\xi|^{-(d/\pp)-m}$, $|\xi|\geq1$,
and $g$ smooth for $|\xi|<1$ shows that 1(i) for $\a=m$
can not be improved. Indeed for this $g$,
$\psi_\pp[{}g{}](t)=\big(|\SD|/(m\pp)\big)^{1/\pp}\cdot t^{-m}$, $t\geq1$,
but $\psi_{\pp,m}[{}g{}](t)\geq 
|\SD|^{1/\pp}\cdot t^{-m}(\log{}t)^{1/\pp}$, $t\geq1$.
\end{remark}
\begin{proof}[Proof of Corollary \ref{CH3_s3.lem.1}]

{\em The case $d\geq2$.\ }The result
follows readily from Corollary~\ref{thm2} 
combined with Theorem~\ref{lem4.1} for $p=\pp=2$, 
$\a=\gamma/2$, $m=1$, $g=\hat{f}$ with $\eps=1/t$. 
Note that the quantitity in the middle in \eqref{eqneww} is 
$(\omega_{2,1,2}[f](\eps))^2$.
The integral in \eqref{eqnewww}
is the true tail integral of the \hbox{F.t.} of $f$, $\psi_2[\hat{f}](t)$.
By Corollary~\ref{thm2}, the exist $C_1,C_2>0$
that depend on $d$ only so that for all $f\in{}L^2(\RD)$
\begin{equation}
\label{eqappr1}
     C_1\,\omega_{2,1,2}[f](1/t)
       \leq \psi_{2,1}[\hat{f}](t)\leq C_2\,\omega_{2,1,2}[f](1/t),\qquad t>0.
\end{equation} 
Set $\alpha=\gamma/2$.
Note that since $0<\gamma<2$
we have $0<\a=\gamma/2<m=1$, and the reult follows
from Theorem~\ref{lem4.1}, cases 1 (i) and (ii), applied to the comparison
function $t^{-\alpha}=t^{-{\gamma/2}}$, $t\geq1$.

{\em The case $d=1$.\ }We prove first that \eqref{eqneww}
implies \eqref{eqnewww}. This follows from a straightforward
modification of the proof of \cite[Lemma 4.2]{BCT}. We prefer to
give the details for the convenience of the reader.
Below, $c$ will denote a positive constant whose precise
value may change from equation to equation and may depend
on $f$ but which is independent of $\eps$ and $t$.
By \eqref{eqParseval} 
$$
\begin{aligned}
              (\omega_{2,1,2}[f](\eps))^2 = \frac4\pi\int_{\RR}
                   \sin^2\bigg(\frac{\eps\xi}{2}\bigg)\,|\hat{f}(\xi)|^2\,d\xi
         \geq c \int_{{\pi}/(2\eps)
                       \leq|\xi|\leq2{\pi}/(2\eps)}
                   |\hat{f}(\xi)|^2\,d\xi.
\end{aligned}
$$
Setting $\tilde{\eps}=2\eps/\pi$ and denoting $\tilde{\eps}$ by $\eps$
again, we find from the upper inequality in \eqref{eqneww}
$$
\begin{aligned}
           c_2\cdot(\pi\eps/2)^\gamma
        \geq   \omega_{2,1,2}[f](\pi\eps/2) 
           \geq c \int_{1/\eps
                       \leq|\xi|\leq2/\eps}
                   |\hat{f}(\xi)|^2\,d\xi
\end{aligned}
$$
which after setting $t=1/\eps$ implies
\begin{equation}\label{eqdyadic}
         \int_{t
                       \leq|\xi|\leq2t}
                   |\hat{f}(\xi)|^2\,d\xi \leq c\cdot t^{-\gamma},\qquad t\geq1.
\end{equation}
Using the latter estimate and 
representing
$$
        \int_{|\xi|\geq t}|\hat{f}(\xi)|^2\,d\xi
          =\sum_{j=0}^\infty 
                 \int_{2^j t\leq|\xi|\leq2^{j+1} t} |\hat{f}(\xi)|^2\,d\xi,
$$
we prove the upper inequality in \eqref{eqnewww}.
In order to prove the lower inquality in \eqref{eqnewww},
we note first that \eqref{eqdyadic} implies
\begin{equation}\label{eq1dyadic}
         \int_{|\xi|\leq r}
                  |\xi|^2\, |\hat{f}(\xi)|^2\,d\xi \leq c_3
                          \cdot r^{2-\gamma},\qquad r\geq1.
\end{equation}
Indeed using the upper inequality in \eqref{eqnewww}
we obtain
$$
\begin{aligned}
      \int_{|\xi|\leq r}
                  |\xi|^2\, |\hat{f}(\xi)|^2\,d\xi 
          &\leq  \int_{|\xi|\leq 1}
                  |\xi|^2\, |\hat{f}(\xi)|^2\,d\xi 
          +\sum_{j=0}^{[\log_2r]+1}
                 \int_{2^j \leq|\xi|\leq2^{j+1} }  |\xi|^2\, |\hat{f}(\xi)|^2\,d\xi\\
         &\leq\const
          +c\sum_{j=0}^{[\log_2r]+1}
                (2^{j+1})^2 \cdot b_2\cdot (2^{j})^{-\gamma}\\
       &\leq c^\prime\cdot r^{2-\gamma}
\end{aligned}
$$
(recall that $\gamma<2$). Next, let $\alpha,\beta>0$ be two
numbers to be chosen later. From the lower inequality
in \eqref{eqneww} and \eqref{eqParseval} using
$|\sin x|\leq x$ and $|\sin x|\leq1$ for $x\geq0$, we obtain
$$
\begin{aligned}
    c_1\cdot \eps^\gamma&\leq
           (\omega_{2,1,2}[f](\eps))^2\\ 
          &= \frac4\pi\int_{\RR}
                   \sin^2\bigg(\frac{\eps\xi}{2}\bigg)\,|\hat{f}(\xi)|^2\,d\xi\\
          &\leq c \bigg(\eps^2\int_{|\xi|\leq\alpha/\eps}
                   |\xi|^2\,|\hat{f}(\xi)|^2\,d\xi
             +\int_{\alpha/\eps\leq|\xi|\leq\beta/\eps}
                   |\hat{f}(\xi)|^2\,d\xi
               +\int_{|\xi|\geq\beta/\eps}
                   |\hat{f}(\xi)|^2\,d\xi\bigg).
\end{aligned}
$$
Using \eqref{eq1dyadic} and the upper inequality in \eqref{eqnewww}
we find
$$
     \int_{\alpha/\eps\leq|\xi|\leq\beta/\eps}
                   |\hat{f}(\xi)|^2\,d\xi 
        \geq c_1c^{-1}\cdot \eps^\gamma 
          -\beta^{-\gamma}b_2\cdot \eps^\gamma
          -\alpha^{\gamma}c_3\cdot \eps^\gamma
$$
which after choosing $\alpha>0$ small enough and $\beta>0$
large enough gives
$$
     \int_{|\xi|\geq\alpha/\eps}
                   |\hat{f}(\xi)|^2\,d\xi 
        \geq\int_{\alpha/\eps\leq|\xi|\leq\beta/\eps}
                   |\hat{f}(\xi)|^2\,d\xi 
        \geq c\cdot\eps^\gamma,\qquad 0<\eps\leq1.
$$
Setting $t=\alpha/\eps$ we prove the lower inequality in \eqref{eqnewww}.

We now derive \eqref{eqneww} from \eqref{eqnewww}.
Using the upper inequality in \eqref{eqnewww}
and employing \eqref{eqtri} and Theorem \ref{lem4.1} 
(recall that in our case $p=p^\prime=2$, $m=1$, $\alpha=\gamma/2\in(0,1)$)
$$
     (\psi_{2}[\hat{f}](t))^2
              \leq (\psi_{2,1}[\hat{f}](t))^2 
               \leq c\cdot t^{-\gamma}.
$$
Combining this with \eqref{eq1} (that holds for $d=1$)
and using \eqref{eq0} we obtain
$$
   (\omega_{2,1,2}[f](1/t))^2 \leq c \cdot(\omega_{2,1,\infty}[f](1/t))^2 
         \leq c^\prime \cdot(\psi_{2,1}[\hat{f}](t))^2
           \leq c^{\prime\prime}\cdot t^{-\gamma}
$$
which proves the upper estimate in \eqref{eqneww}.
It remains to prove the lower estimate in \eqref{eqneww}.
 Note that the two-sided estimate \eqref{eqnewww}
implies that for $A>1$ large enough 
$$
         \bigg( \int_{|\xi|\geq t} -
                     \int_{|\xi|\geq At} \bigg)\,|\hat{f}(\xi)|^2\,d\xi
          =  \int_{t\leq|\xi|\leq At} |\hat{f}(\xi)|^2\,d\xi \geq c\cdot t^{-\gamma}
$$
or after setting $\eps=\pi/(2t)$
\begin{equation}\label{eq2dyadic}
   \int_{\pi/(2\eps)\leq|\xi|\leq A\pi/(2\eps)}
               |\hat{f}(\xi)|^2\,d\xi \geq c\cdot \eps^{\gamma},
                     \qquad 0<\eps\leq1.
\end{equation}
On the other hand again using \eqref{eqParseval} 
$$
\begin{aligned}
              (\omega_{2,1,2}[f](\eps))^2 = \frac4\pi\int_{\RR}
                   \sin^2\bigg(\frac{\eps\xi}{2}\bigg)\,|\hat{f}(\xi)|^2\,d\xi
         \geq c \int_{{\pi}/(2\eps)
                       \leq|\xi|\leq3{\pi}/(2\eps)}
                   |\hat{f}(\xi)|^2\,d\xi
\end{aligned}
$$
and replacing $\eps$ with $\eps/3^j$, $j\in\NN$, we find
$$
              (\omega_{2,1,2}[f](3^{-j}\eps))^2 
           \geq c \int_{3^j{\pi}/(2\eps)
                       \leq|\xi|\leq3^{j+1}{\pi}/(2\eps)}
                   |\hat{f}(\xi)|^2\,d\xi,\qquad j=0,1,2,\cdots.
$$
Choosing
$N=[\log_3 A]$ we then obtain
$$
\begin{aligned}
          \sum_{j=0}^N (\omega_{2,1,2}[f](3^{-j}\eps))^2 
                 &\geq c\cdot
           \sum_{j=0}^N  \int_{3^j{\pi}/(2\eps)
                       \leq|\xi|\leq3^{j+1}{\pi}/(2\eps)}
                   |\hat{f}(\xi)|^2\,d\xi \\
       & \geq c\cdot
           \sum_{j=0}^N  \int_{{\pi}/(2\eps)
                       \leq|\xi|\leq A{\pi}/(2\eps)}
                   |\hat{f}(\xi)|^2\,d\xi \\
        &\geq c^\prime\cdot \eps^\gamma,\qquad 0<\eps\leq1,
\end{aligned}
$$
where we have used \eqref{eq2dyadic}. Therefore
$$
      \liminf_{\eps\to0}\eps^{-\gamma}
              \sum_{j=0}^N \big(\omega_{2,1,2}[f](3^{-j}\eps)\big)^2\geq c>0
$$
and hence as $N\in\NN$ is fixed, there exists at least one $J\in\{1,\cdots,N\}$
such that 
$$
   \tilde{c}:=
       \liminf_{\eps\to0} \eps^{-\gamma}
                   \big(\omega_{2,1,2}[f](3^{-J}\eps)\big)^2 > 0.
$$
Then 
$$
       \big(\omega_{2,1,2}[f](3^{-J}\eps)\big)^2\geq \tilde{c}\cdot \eps^\gamma,
                 \qquad 0<\eps\leq1,
$$
and denoting $3^{-J}\eps$ by $\eps$ we finish the proof of the lower 
inequality in \eqref{eqneww}.
\end{proof}
We note finally that the proof of  Corollary \ref{CH3_s3.lem.1}
for $d=1$ 
can be modified 
to give an alternative proof of  Corollary \ref{CH3_s3.lem.1}
 also for all $d\geq2$
from scratch (in this connection, see an explanation of an argument
from \cite{BCT} given in the proof of \cite[Lemma~3.4.1]{thesis}).
\begin{remark}
\label{remthf}
After the above general discussion
it is not difficult to understand why 
Corollary~\ref{CH3_s3.lem.1}
fails for $\gamma=2$.
Recall that in \eqref{eqappr1}, $p=\pp=2$ and $m=1$.
Let $\hat{f}\in{}L^2(\RD)$ be
defined by $\hat{f}(\xi):=|\xi|^{-(d/2)-1}$ for $|\xi|\geq1$
and smooth for $|\xi|<1$. 
Let $f\in{}L^2(\RD)$ be the inverse \hbox{F.t.} of this $\hat{f}$.
It is easy to check that for some $\tilde{b}_1,\tilde{b}_2>0$ that depend
on $\hat{f}$
\begin{equation}
\label{eqtss}
     \tilde{b}_1\,t^{-1}(\log{t})^{1/2} \leq\psi_{2,1}[\hat{f}](t)
           \leq \tilde{b}_2\,t^{-1}(\log{t})^{1/2},\qquad t\geq2
\end{equation}
whereas for certain $\tilde{c}_1,\tilde{c}_2>0$ that depend
on $\hat{f}$
$$
     \tilde{c}_1\,t^{-1} \leq\psi_{2}[\hat{f}](t)
           \leq \tilde{c}_2\,t^{-1},\qquad t\geq2.
$$
Note that by \eqref{eqappr1} and \eqref{eqtss}
for $c_1,c_2>0$ that depend on $f$
$$
     c_1\,t^{-1}(\log{t})^{1/2} \leq\omega_{2,1,2}[f](1/t)
           \leq c_2\,t^{-1}(\log{t})^{1/2},\qquad t\geq2.
$$
This shows that Corollary~\ref{CH3_s3.lem.1}
fails for $\gamma=2$.
It is only true in the case $\gamma=2$ that the upper
estimate in \eqref{eqneww} implies 
the upper estimate in \eqref{eqnewww}
(simply because for any $f\in{}L^2(\RD)$ and all $t>0$,
$\psi_2[\hat{f}](t)\leq{}\psi_{2,1}[\hat{f}](t)\leq{}c(2,d,1)\,\omega_{2,1,2}[f](1/t)$
in view of \eqref{eqtri} and \eqref{eq2p}).
The fact that the lower estimate in \eqref{eqnewww} need not hold
is shown by considering the
example of $\hat{f}\in{}C_0^\infty(\RD)$, $\hat{f}\not\equiv0$ (for
which $\psi_2[\hat{f}](t)=0$ identically for large $t$).
Finally, the first example of this remark 
shows the upper estimate in \eqref{eqneww} for $\gamma=2$
need not follow even
from a two-sided estimate in \eqref{eqnewww}.
\end{remark}
\bibliographystyle{amsalpha}

\end{document}